\documentclass[12pt]{amsart}
\newtheorem{theorem}{Theorem}[section]
\newtheorem{lemma}[theorem]{Lemma}
\newtheorem{proposition}[theorem]{Proposition}
\newtheorem{definition}{Definition}[section]
\newtheorem{corollary}[theorem]{Corollary}

\newcommand{\conv}{\mbox{{\rm conv}}}
\newcommand{\fvec}{$\{f_\emptyset^\diamond,f_0^\diamond,f_1^\diamond,\ldots,f_{d-1}^\diamond\}$}

\begin{document}

\title{Flag vectors of multiplicial polytopes}

\author{Margaret M. Bayer}
\address{\hskip-\parindent
        Department of Mathematics\\
        University of Kansas\\
        Lawrence, KS  66045-7523}
\email{bayer@math.ukans.edu}

\thanks{Research at MSRI is supported in part by NSF grant DMS-9810361.
This research was also supported by a University of Kansas General Research
Fund Grant.}

\begin{abstract}
Bisztriczky introduced the {\em multiplex} as a generalization of the simplex.
A polytope is multiplicial if all its faces are multiplexes.
In this paper it is proved that the flag vectors of multiplicial polytopes
depend only on their face vectors.  A special class of multiplicial
polytopes, also discovered by Bisztriczky, is comprised of the {\em ordinary
polytopes}.  These are a natural generalization of the cyclic polytopes.
The flag vectors of ordinary polytopes are determined.  This is used to give
a surprisingly simple formula for the $h$-vector of the ordinary $d$-polytope
with $n+1$ vertices and characteristic $k$:
$h_i=\binom{k-d+i}{i}+(n-k)\binom{k-d+i-1}{i-1}$, for $i\le d/2$.
In addition, a construction is given for 4-dimensional multiplicial polytopes
having two-thirds of their vertices on a single facet, answering a question
of Bisztriczky.
\end{abstract}

\maketitle

\section{Introduction}
Convex polytopes arise in many areas of mathematics and its applications.
Of great interest are 
their combinatorial aspects, and, in particular, the numbers of faces
of convex polytopes.
In this area, the biggest result was the characterization of the $f$-vector
(the sequence of face
numbers) of simplicial convex polytopes in 1980 by Billera and Lee 
\cite{billera-lee} and Stanley \cite{stanley}.
Simplicial polytopes are those for which every proper face is a simplex.
For nonsimplicial polytopes we have seen no such success.
Contributing to our failure is the lack of constructions of polytopes 
with sufficiently varied facial structures.

In \cite{Bisz-mult} Bisztriczky introduced the {\em multiplex} as a 
generalization of the simplex.
There is a $d$-dimensional multiplex with $v$ vertices for each $v\ge d+1$.
Multiplexes share many of the nice properties of simplices; in particular
they are self-dual and every face and every quotient of a multiplex is a
multiplex.
It is natural, therefore, to study the combinatorics of {\em multiplicial}
polytopes, polytopes whose faces are multiplexes.
All polygons are multiplexes, so all 3-dimensional polytopes are multiplicial.
The story begins, therefore, in dimension 4.
This is also where the open problems on face numbers of general polytopes
begin.

In \cite{Bisz} Bisztriczky introduced a fascinating class of multiplicial
polytopes, called {\em ordinary} polytopes.
Multiplexes and ordinary polytopes have been further studied in \cite{Dinh}
and \cite{Bayer-Bru-Ste}.
This paper continues the study of $f$-vectors and flag vectors of 
multiplicial polytopes.
In addition, we give a construction of multiplicial 4-polytopes with large
facets.

We turn now to the necessary definitions and background theorems.

\begin{definition}[\cite{Bisz-mult}]  {\em
A $d$-dimensional {\em multiplex} is a polytope with an ordered list of 
vertices,
$x_0$, $x_1$, \ldots, $x_n$, with facets $F_0$, $F_1$, \ldots, $F_n$ given by
$$F_i = \conv\{x_{i-d+1}, x_{i-d+2}, \ldots, x_{i-1}, x_{i+1}, x_{i+2}, \ldots,
x_{i+d-1}\},$$
with the conventions that $x_i=x_0$ if $i<0$, and $x_i=x_n$ if $i>n$.
}\end{definition}

\begin{theorem}
[\cite{Bisz-mult}]
\label{multi-backgr}
$~$
\begin{enumerate}
\item For every $d$ and $n$ with $n\ge d\ge 2$,
      there exists a $d$-dimensional multiplex $M^{d,n}$ with $n+1$
      vertices. 
\item Every multiplex is self-dual.
\item Every face and every quotient of a multiplex is a multiplex.
\item The number of $i$-dimensional faces of $M^{d,n}$
      is $\binom{d+1}{i+1}+(n-d)\binom{d-1}{i}$.
\end{enumerate}
\end{theorem}

By definition the combinatorial type of $M^{d,n}$ is completely determined
by $d$ and $n$.  If $d=n$, of course, the multiplex is the simplex.

In the study of nonsimplicial polytopes, the $f$-vector does not carry
enough combinatorial information.
We are interested in incidences of faces as well.
A chain of faces, $\emptyset\subset F_1\subset F_2\subset \cdots
\subset F_r\subset P$ is an {\em $S$-flag} if
$S=\{\dim F_1, \dim F_2, \ldots, \dim F_r\}$.
The number of $S$-flags of a polytope $P$ is written $f_S(P)$, and
the length $2^d$ vector $(f_S(P))_{S\subseteq\{0, 1, \ldots, d-1\}}$
is the {\em flag vector} of $P$.
The flag vector restricts to the $f$-vector by considering only the
singleton sets $S$.
The flag vector of the multiplex is a familiar object.

\begin{theorem} [\cite{Bayer-Bru-Ste}]
\label{flag-multiplex}
The multiplex $M^{d,n}$ has the same flag vector as the $(d-2)$-fold
pyramid over the $(n-d+3)$-gon.
The common flag vector is given by
\begin{eqnarray}
\label{multiplex-flag}
f_S&=&\binom{d+1}{\mbox{$s_1+1$, $s_2-s_1$, \ldots, $s_r-s_{r-1}$,
  $d-s_r$}}\\
  & & \times\left[ 1+\frac{n-d}{(d+1)d(d-1)}\sum_{j=1}^r
  (s_j+1)(s_{j+1}-s_j)(s_{j+1}-1) \right], \nonumber
\end{eqnarray}
where $S=\{s_1,s_2,\ldots, s_r\}$, $s_1<s_2<\cdots<s_r$, and $s_{r+1}=d$.
\end{theorem}

\section{Results for general multiplicial polytopes}
There are two natural classes of polytopes that could be called 
``multiplicial''.  
One might call a polytope multiplicial if all of its proper faces are
multiplexes.
However, a  multiplex is really a polytope endowed with a special 
ordering of its vertices.
It is natural to take the ordering into account when defining multiplicial,
and Bisztriczky indeed does this.
We will use both concepts in this paper.
\begin{definition}
{\em
A polytope $P$ is {\em (weakly) multiplicial} if and only if every proper 
face of $P$ is a multiplex.
}
\end{definition}
\begin{definition}
{\em
A polytope $P$ is {\em order-multiplicial} if and only if for some ordering
$v_0$, $v_1$, \ldots, $v_n$ of the vertices of $P$, every proper 
face of $P$ is a multiplex with the induced ordering of its vertices.
}
\end{definition}

Let ${\mathcal P}^d$ be the set of $d$-polytopes in ${\bf R}^d$.
Let $f_S:{\mathcal P}^d\rightarrow {\bf N}$ be the function that gives the
number of $S$-flags of $P$ for a $d$-polytope $P$.
In particular $f_i$ is the function that gives the number of $i$-faces of $P$.
Write $f_S^\diamond$ for the restriction of the $S$-flag function to 
multiplicial polytopes.
For a multiplicial polytope, the flag vector of each proper face depends
only on its number of vertices.
This leads to the following result.

\begin{theorem}\label{flag-multiplicial-1}
Fix $d\ge 2$, and let $S$ be a nonempty subset of $\{0,1,\ldots, d-1\}$ with 
maximum element $t$.
There exist rational numbers $a$ and $b$ such that 
$f_S^\diamond=af_t^\diamond+bf_{0,t}^\diamond$.

For $S=\{s_1,s_2,\ldots, s_r\}$, with $s_1<s_2<\cdots<s_r=t$,
\begin{eqnarray*}
a&=&\binom{s_r+1}{\mbox{$s_1+1$, $s_2-s_1$, \ldots, $s_r-s_{r-1}$}}\\
  & & \times\left[ 1-\frac{1}{s_r(s_r-1)}\sum_{j=1}^{r-1}
  (s_j+1)(s_{j+1}-s_j)(s_{j+1}-1) \right],
\end{eqnarray*}
and 
\begin{eqnarray*}
b&=&\binom{s_r+1}{\mbox{$s_1+1$, $s_2-s_1$, \ldots, $s_r-s_{r-1}$ }}\\
  & & \times\left[\frac{1}{(s_r+1)s_r(s_r-1)}\sum_{j=1}^{r-1}
  (s_j+1)(s_{j+1}-s_j)(s_{j+1}-1) \right].
\end{eqnarray*}
\end{theorem}
\begin{proof}
We first observe that the statement holds trivially for $|S|=1$.
For $|S|\ge 2$, we use the formula for the flag vectors of multiplexes
(see Theorem~\ref{flag-multiplex}).
Suppose $S=\{s_1,s_2,\ldots, s_r\}$, with $s_1<s_2<\cdots<s_r$, and 
write $s_{r+1}=d$.
For the $d$-multiplex $M^{d,n}$ with $n+1$ vertices, 
$f_S(M^{d,n})=a_{S,d}+b_{S,d}(n+1)=a_{S,d}+b_{S,d}f_0(M^{d,n})$,
where
\begin{eqnarray*}
a_{S,d}&=&\binom{d+1}{\mbox{$s_1+1$, $s_2-s_1$, \ldots, $s_r-s_{r-1}$,
  $d-s_r$}}\\
  & & \times\left[ 1-\frac{1}{d(d-1)}\sum_{j=1}^r
  (s_j+1)(s_{j+1}-s_j)(s_{j+1}-1) \right],
\end{eqnarray*}
and 
\begin{eqnarray*}
b_{S,d}&=&\binom{d+1}{\mbox{$s_1+1$, $s_2-s_1$, \ldots, $s_r-s_{r-1}$,
  $d-s_r$}}\\
  & & \times\left[\frac{1}{(d+1)d(d-1)}\sum_{j=1}^r
  (s_j+1)(s_{j+1}-s_j)(s_{j+1}-1) \right].
\end{eqnarray*}
Now fix $d\ge 2$, and let $S$ be a subset of $\{0,1,\ldots, d-1\}$ with 
$|S|\ge 2$ and maximum element $t$.
Let $P$ be a multiplicial $d$-polytope.
Every $t$-face $F$ of $P$ is a $t$-dimensional multiplex, so
$\displaystyle
f_{S\setminus\{t\}}(F)=a_{S\setminus\{t\},t}+b_{S\setminus\{t\},t}f_0(F)$.
Then 
\begin{eqnarray*}
f_S(P)&=& \sum_{\mbox{{\scriptsize $F$ $t$-face of $P$}}} f_{S\setminus\{t\}}
(F)\\
      &=& \sum_{\mbox{{\scriptsize $F$ $t$-face of $P$}}} a_{S\setminus\{t\},t}
          +b_{S\setminus\{t\},t}f_0(F)\\
      &=& a_{S\setminus\{t\},t}f_t(P) +b_{S\setminus\{t\},t}f_{0,t}(P).\\
\end{eqnarray*}

\vspace*{-24pt}

\end{proof}

\begin{theorem}\label{flag-multiplicial-2}
For each $S\subseteq\{0,1,\ldots,d-1\}$, the function $f_S^\diamond$ 
is a linear combination of the constant function $f_\emptyset^\diamond(P)=1$ 
and the face number functions $f_i^\diamond$, $0\le i\le d-2$.
\end{theorem}

\begin{proof}
We have already shown that every $f_S^\diamond$ is in the linear span of \break
$\{f_\emptyset^\diamond, f_0^\diamond,f_1^\diamond, \ldots, f_{d-1}^\diamond,
f_{0,1}^\diamond,f_{0,2}^\diamond,\ldots, f_{0,d-1}^\diamond\}$.
Now we show that every $f_{0,t}^\diamond$ is in the linear span of 
$\{f_\emptyset^\diamond, f_0^\diamond,f_1^\diamond, \ldots, f_{d-1}^\diamond\}$.
Of course, Euler's formula enables us to drop $f_{d-1}^\diamond$.

The proof is by downward induction on $t$.
We check first that it holds for $t=d-1$.
For a multiplicial $d$-polytope $P$ and any $(d-1)$-face $F$ of $P$, 
$f_0(F)=f_{d-2}(F)$, since $F$ is a multiplex.
Then
\begin{eqnarray*}
f_{0,d-1}(P)&=&\sum_{\mbox{{\scriptsize $F$ $(d-1)$-face of $P$}}} f_0(F)\\
&=&\sum_{\mbox{{\scriptsize $F$ $(d-1)$-face of $P$}}} f_{d-2}(F)=f_{d-2,d-1}(P)=2f_{d-2}(P).
\end{eqnarray*}

Now assume for all $k>t$, $f_{0,k}^\diamond$ is in the span of \fvec.
Since all $t$-faces of a multiplicial polytope are multiplexes, we get (as in
the $t=d-1$ case) $f_{0,t}^\diamond=f_{t-1,t}^\diamond$.
Now the generalized Dehn-Sommerville equations say that
$f_{t-1,t}-f_{t-1,t+1}+\cdots +(-1)^{d-1-t}f_{t-1,d-1}=(1-(-1)^{d-t})f_{t-1}$.
(This results from applying Euler's formula to the quotient of the polytope $P$
by each $(t-1)$-face.)
By Theorem~\ref{flag-multiplicial-1} each $f_{t-1,k}^\diamond$ is in the span of
$\{f_0^\diamond,f_{0,k}^\diamond\}$, so
by the induction hypothesis, $f_{t-1,k}^\diamond$ is in the span of \fvec\ for 
all $k\ge t+1$, so $f_{t-1,t}^\diamond=f_{0,t}^\diamond$ is also in this linear
span.
\end{proof}

\begin{corollary}\label{dim-multiplicial}
The linear span of the flag vectors of all multiplicial $d$-polytopes 
is of dimension at most $d$.
\end{corollary}

Computation shows that the dimension is in fact $d$ for $d\le 37$, with the
{\em ordinary polytopes} (see next section) providing a set of spanning flag
vectors.
Note that the induction argument in the proof of the theorem does not
yield an attractive formula for the flag numbers in terms of the face
numbers.

\section{Flag vectors of ordinary polytopes}
Of special interest among order-multiplicial polytopes are the 
{\em ordinary polytopes}, discovered by Bisztriczky \cite{Bisz}.
Given an ordered set $V=\{x_0,x_1,\ldots,x_n\}$, a subset $Y\subseteq V$ is
called a {\em Gale subset} if between any two elements of $V\setminus Y$
there is an even number of elements of $Y$.
A polytope $P$ with ordered vertex set $V$ is a {\em Gale polytope}
if the set of vertices of each facet is a Gale subset.
\begin{definition}[\cite{Bisz}] {\em
An {\em ordinary polytope} is a Gale polytope such that each facet is a
multiplex with the induced order on the vertices.
}\end{definition}
For background on ordinary polytopes, see \cite{Bayer-Bru-Ste,Bisz,Dinh}.
Three-dimensional ordinary polytopes are quite different from those in higher
dimensions, so we will exclude them from consideration.
Higher-dimensional ordinary polytopes provide a natural combinatorial
generalization of cyclic polytopes, which have played such an important role
in the combinatorial study of simplicial polytopes.
\begin{theorem}[\cite{Bisz}]
Let $P$ be an ordinary $d$-polytope with ordered vertices
$x_0$, $x_1$, \ldots, $x_n$.
Assume $n\ge d\ge 4$.
\begin{enumerate}
\item If $d$ is even, then $P$ is cyclic.  
\item If $d$ is odd, then there exists an integer $k$ ($d\le k\le n$) such that
      the vertices sharing an edge with $x_0$ are exactly $x_1$, $x_2$, \ldots,
      $x_k$, and the vertices sharing an edge with $x_n$ are exactly
      $x_{n-1}$, $x_{n-2}$, \ldots, $x_{n-k}$.
\end{enumerate}
\end{theorem}

From now on we will restrict our attention to ordinary polytopes of odd
dimension, $d=2m+1\ge 5$.
The integer $k$ guaranteed by this theorem is called the {\em characteristic}
of the ordinary polytope.
The three parameters $d$, $k$, and $n$ completely determine the combinatorial
structure of the ordinary polytope, which is denoted $P^{d,k,n}$.
When $k=d$, the ordinary polytope is a multiplex.
When $k=n$, the ordinary polytope is a cyclic polytope.

Dinh (\cite{Dinh}) gives an inductive geometric
construction of $P^{d,k,n}$ starting with the cyclic polytope $P^{d,k,k}$.
He derives from this construction formulas for the face numbers $f_i$ in
terms of $d$, $k$, and $n$.
We extend this method to obtain formulas for
the flag numbers $f_{0,t}$ similar to Dinh's $f$-vector formulas.

Dinh successively introduces a new vertex $x_n$ to produce the ordinary
polytope $P^{d,k,n}$ as $\conv(P^{d,k,n-1}\cup\{x_n\})$.
He identifies and counts two classes of faces.
A {\em new $i$-face of $P^{d,k,n}$} is an $i$-dimensional face of 
$P^{d,k,n}$ that does not contain an $i$-face of $P^{d,k,n-1}$.
An {$i$-face of $P^{d,k,n-1}$ destroyed by $x_n$} is an $i$-face of 
$P^{d,k,n-1}$ not contained in an $i$-face of $P^{d,k,n}$.
The other faces of $P^{d,k,n-1}$ are preserved as faces in $P^{d,k,n}$,
either intact or with the vertex $x_n$ added in the affine span.

\begin{proposition}[\cite{Dinh}]
The new $i$-faces of $P^{d,k,n}$ are all simplices.
The number of them is
$$ \alpha_i(d,k)= \sum_{j=i-m}^{\lfloor i/2\rfloor} 2N(k-1,j,i)-N(k-2,j,i),$$
where 
$$N(s,t,u)=\binom{u-t}{t}\binom{s-u+t}{u-t}+\binom{u-1-t}{t}\binom{s-u+t}
  {u-1-t}.$$
When $i\le m$ this simplifies to
$ \alpha_i(d,k)= \binom{k-1}{i}+\binom{k-2}{i-1} $.
\end{proposition}

\begin{theorem}
The new $i$-faces of $P^{d,k,n}$ 
and the destroyed $i$-faces of $P^{d,k,n-1}$ are all simplices.
For $0\le i\le 2m=d-1$, the number of new $i$-faces minus the number of 
destroyed $i$-faces, $c_i(d,k)=f_i(P^{d,k,n})-f_i(P^{d,k,n-1})$, is given by
$$c_i(d,k)= \sum_{j=0}^{m-1}\binom{m+1}{i-j}\binom{k-m-2}{j}.$$
When $i\le m-1$ this simplifies to
$ c_i(d,k)= \binom{k-1}{i}$.
\end{theorem}

Dinh \cite{Dinh} first computed the quantities $c_i(d,k)$, but he gave a more 
complicated expression.
His description of the new and destroyed faces, however, lends itself to
the technique used to count the faces of cyclic polytopes in \cite{McMullen}.
After some manipulation, this yields the relatively nice form stated here.
It would be good to have a bijective proof of the formula.

We also need the $f$-vectors of the cyclic polytopes.
(See \cite{McMullen}.)
\begin{proposition}
Let $m=\lfloor d/2\rfloor$.
The cyclic $d$-polytope with $k+1$ vertices has $f$-vector given by
\begin{eqnarray*}
f_i=\phi_i(d,k)&=& \sum_{j=0}^m \left[ \binom{j}{d-1-i}+\binom{d-j}{d-1-i}\right]
       \binom{k-d+j}{j}\\
       & & -\chi(\mbox{$d$ even})\binom{m }{2m-1-i}\binom{k-m}{m}
\end{eqnarray*}
When $i\le m-1$ this simplifies to
$ \phi_i(d,k)= \binom{k+1}{i+1}$.
\end{proposition}

Here are the $f$-vectors of the ordinary polytopes, as computed by Dinh.
\begin{theorem}[\cite{Dinh}]
Let $n\ge k\ge d=2m+1\ge 5$.
The number of $i$-dimensional faces of the ordinary polytope $P^{d,k,n}$ is
$$f_i(P^{d,k,n})=\phi_i(d,k)+(n-k)c_i(d,k),$$
\end{theorem}

Here are the analogous formulas for the flag numbers $f_{0,i}$.

\begin{theorem} For $ n\ge k\ge d=2m+1\ge 5$ and $1\le i\le d-1$,
\begin{eqnarray*}
f_{0,i}(P^{d,k,n})&=&(i+1)\phi_i(d,k)\\ &+& (n-k)[(i+1)c_i(d,k)
                     +\phi_{i-1}(d-1,k-1)-\alpha_i(d,k)].
\end{eqnarray*}
\end{theorem}
\begin{proof}
To get the number of incidences of vertices and $i$-faces for $P^{d,k,n}$,
start with $f_{0,i}(P^{d,k,n-1})$, add the number of vertices on each new
$i$-face, subtract the number of vertices on each destroyed $i$-face,
and add one for each preserved $i$-face that now contains the new vertex 
$x_n$.
All the new and destroyed faces are simplices, and the net gain in $i$-faces
is $c_i$, so the contribution to the change in $f_{0,i}$ from the new
and destroyed faces is $(i+1)c_i$.
It remains to consider the preserved $i$-faces that now contain vertex $x_n$.

\begin{lemma} The number of $i$-faces of $P^{d,k,n}$ containing the vertex
$x_n$ is \break
$\phi_{i-1}(d-1,k-1)$, the number of $(i-1)$-faces of the cyclic
$(d-1)$-polytope with $k$ vertices.
\end{lemma}
\noindent {\em Proof of lemma:}
From the definition of ordinary polytope, it is clear that
reversing the ordering of the vertices produces the same polytope.
So the number of $i$-faces of $P^{d,k,n}$ containing the vertex $x_n$ is
the same as the number of $i$-faces of $P^{d,k,n}$ containing the vertex $x_0$.
In Dinh's inductive construction of $P^{d,k,n}$ from $P^{d,k,n-1}$, none of 
the new $i$-faces contain $x_0$, and none of the destroyed $i$-faces contain
$x_0$.
So the number of $i$-faces of $P^{d,k,n}$ containing the vertex $x_0$ is the 
same as the number of $i$-faces of the cyclic polytope $P^{d,k,k}$ containing 
the vertex $x_0$.
The $i$-faces of $P^{d,k,k}$ containing $x_0$ are in one-to-one correspondence
with the $(i-1)$-faces of the link of vertex $x_0$ in $P^{d,k,k}$.
By Gale's evenness condition for the facets of the cyclic polytope, the link
of vertex $x_0$ in the cyclic $d$-polytope with $k+1$ vertices is the cyclic
$(d-1)$-polytope with $k$ vertices.
\ \hfill{\em End of proof of lemma}

\vspace*{6pt}

All the new $i$-faces of $P^{d,k,n}$ contain the vertex $x_n$, so the number
of preserved $i$-faces that now contain vertex $x_n$ is
$\phi_{i-1}(d-1,k-1)-\alpha_i$.
Since the cyclic polytope $P^{d,k,k}$ is simplicial, 
$f_{0,i}(P^{d,k,k})=(i+1)f_i(P^{d,k,k})=(i+1)\phi_i(d,k)$.  Thus
\begin{eqnarray*}
f_{0,i}(P^{d,k,n})&=& f_{0,i}(P^{d,k,n-1}) +(i+1)c_i
                     +\phi_{i-1}(d-1,k-1)-\alpha_i\\
                 &=& f_{0,i}(P^{d,k,k})+(n-k)[(i+1)c_i
                     +\phi_{i-1}(d-1,k-1)-\alpha_i]\\
                 &=& (i+1)\phi_i(d,k)+(n-k)[(i+1)c_i
                     +\phi_{i-1}(d-1,k-1)-\alpha_i].
\end{eqnarray*}

\vspace*{-12pt}

\end{proof}

The numbers $\alpha_i$, $c_i$ and $\phi_i$ are independent of $n$.
So for fixed $d$ and $k$, the flag vectors of the polytopes $P^{d,k,n}$
lie on a line.

\section{Toric $h$-vectors}

We turn now to the toric $h$-vectors of ordinary polytopes.
In the characterization of $f$-vectors of simplicial polytopes
\cite{billera-lee,stanley}, a crucial role was played by the $h$-vector,
the image under a certain linear map of the $f$-vector.
The $h$-vector can be interpreted as the sequence of homology ranks of the toric
variety associated to the polytopes.
In this sense (using intersection homology) the $h$-vector is also defined for 
nonsimplicial polytopes, but it depends on the flag vector, not just the 
$f$-vector.
(See \cite{bay-ehr} for formulas in terms of the flag vector.)
Following Stanley \cite{sta87} we define the toric $h$-vector (and $g$-vector)
of any polytope (or Eulerian poset),
first encoding the $h$-vector and $g$-vector as polynomials:
$h(P,t)=\sum_{i=0}^d h_it^{d-i}$ and
$g(P,t)=\sum_{i=0}^{\lfloor d/2\rfloor}g_it^i$, with the relations $g_0=h_0$
and $g_i=h_i-h_{i-1}$ for $1\le i\le d/2$.
Then the $h$-vector and $g$-vector are defined by the recursion
\begin{enumerate} 
\item $g(\emptyset,t)=h(\emptyset,t)=1$, and
\item $\displaystyle h(P,t)=\sum_{\genfrac{}{}{0pt}{}{\mbox{{\scriptsize $G$ face of $P$}}}
       {\mbox{{\scriptsize $G\not=P$}}}} g(G,t)(t-1)^{d-1-\dim G}$.
\end{enumerate}

The toric $h$-vectors of ordinary polytopes have a very simple form.

\pagebreak

\begin{theorem} For $ n\ge k\ge d=2m+1\ge 5$ and $1\le i\le m$,
$$h_i(P^{d,k,n})=\binom{k-d+i}{i}+(n-k)\binom{k-d+i-1}{i-1}.$$
\end{theorem}
\begin{proof}
It follows from Theorem~\ref{flag-multiplex} that the $g$-polynomial of
the multiplex $M^{d,n}$ is the same as the $g$-polynomial of the 
$(d-2)$-fold pyramid over the $(n-d+3)$-gon.
This is the same as the $g$-polynomial of the $(n-d+3)$-gon itself, so 
$g(M^{d,n},t)=1+(n-d)t$.
So for any multiplicial $d$-polytope $P$ with $n+1$ vertices,
\begin{eqnarray*}
h(P,t)&=&\sum_{\genfrac{}{}{0pt}{}{\mbox{{\scriptsize $G$ face of $P$}}}
       {\mbox{{\scriptsize $G\not=P$}}}} g(G,t)(t-1)^{d-1-\dim G}\\
       &=& (t-1)^d+(n+1)(t-1)^{d-1} \\
       & &\mbox{}+\sum_{i=1}^{d-1}
        \sum_{\genfrac{}{}{0pt}{}{\mbox{{\scriptsize $G$ face of $P$}}}
        {\mbox{{\scriptsize $\dim G=i$}}}}
       (1+(f_0(G)-1-i)t)(t-1)^{d-1-i}\\
   &=& (t-1)^d+(n+1)(t-1)^{d-1} +\sum_{i=1}^{d-1}f_i(P)(t-1)^{d-1-i}(1-(i+1)t)\\
      & &\mbox{} +\sum_{i=1}^{d-1}
       \sum_{\genfrac{}{}{0pt}{}{\mbox{{\scriptsize $G$ face of $P$}}}
         {\mbox{{\scriptsize $\dim G=i$}}}} f_0(G)\,t(t-1)^{d-1-i}.
\end{eqnarray*}
So for any multiplicial $d$-polytope $P$ with $n+1$ vertices,
\begin{eqnarray*}
 h(P,t)  &=& (t-1)^d+(n+1)(t-1)^{d-1} +\sum_{i=1}^{d-1}f_i(P)(t-1)^{d-1-i}\\
   & &\mbox{} +\sum_{i=1}^{d-1}\left(f_{0i}(P)-(i+1)f_i(P)\right)t(t-1)^{d-1-i}.
\end{eqnarray*}
Now for $P$ an ordinary polytope, substitute 
\begin{eqnarray*}
n+1 &=& (k+1)+(n-k)\\
f_i(P) &=& \phi_i(d,k)+(n-k)c_i(d,k)\\
f_{0i}(P)-(i+1)f_i(P)&=&(n-k)[\phi_{i-1}(d-1,k-1)-\alpha_i(d,k)]
\end{eqnarray*}
The result is
\begin{eqnarray*}
h(P^{d,k,n},t)&=&
      (t-1)^d+(k+1)(t-1)^{d-1} +\sum_{i=1}^{d-1}\phi_i(d,k)(t-1)^{d-1-i}\\
      & &\mbox{} + (n-k)\left[(t-1)^{d-1} +\sum_{i=1}^{d-1}c_i(d,k)(t-1)^{d-1-i}
          \right.\\
      & &\left.\mbox{} +\sum_{i=1}^{d-1}[\phi_{i-1}(d-1,k-1)-
	   \alpha_i(d,k)]t(t-1)^{d-1-i}\right]\\
      &=& h(P^{d,k,k})+
       (n-k)\left[\sum_{i=0}^{d-1}c_i(d,k)(t-1)^{d-1-i}\right.\\
      & &\left.\mbox{} +\sum_{i=1}^{d-1}[\phi_{i-1}(d-1,k-1)-
	   \alpha_i(d,k)]t(t-1)^{d-1-i}\right].\\
\end{eqnarray*}
The functions $c_i(d,k)$, $\phi_{i-1}(d-1,k-1)$, and $\alpha_i(d,k)$ have simple
expressions when $i$ is small relative to $d$.
This enables us to compute the high degree terms of $h(P^{d,k,n})$.
Let $A=(h(P^{d,k,n})-h(P^{d,k,k}))/(n-k)$, that is,
$$A= \sum_{i=0}^{d-1}c_i(d,k)(t-1)^{d-1-i}
      +\sum_{i=1}^{d-1}[\phi_{i-1}(d-1,k-1)-
	   \alpha_i(d,k)]t(t-1)^{d-1-i}.$$
Let $A_j$ be the coefficient of $t^j$ in $A$.
To compute $A_j$ for $j\ge m+1$, note that $\deg(t-1)^{d-1-i}\ge j$ if and only
if $i\le d-1-j$, and $\deg t(t-1)^{d-1-i}\ge j$ if and only if $i\le d-j$.
Thus 
\begin{eqnarray*}
A_j &=& \sum_{i=0}^{d-1-j}(-1)^{d-1-i-j}\binom{d-1-i}{j}c_i(d,k)\\
    & & \mbox{} + \sum_{i=1}^{d-j}(-1)^{d-i-j}\binom{d-1-i}{j-1}
	 [\phi_{i-1}(d-1,k-1)-\alpha_i(d,k)].
\end{eqnarray*}
For $j\ge m+1$ and $i\le d-1-j$, $i\le d-m-2=m-1$, so we need $c_i(d,k)$ only 
for $i\le m-1$.
For $j\ge m+1$ and $i\le d-j$, $i\le d-m-1=m$, so we need 
$\phi_{i-1}(d-1,k-1)-\alpha_i(d,k)$ only for $i\le m$.
For $i\le m-1$, $c_i(d,k)=\binom{k-1}{i}$.
For $1\le i\le m$, $\phi_{i-1}(d-1,k-1)=\binom{k}{i}$ and
$\alpha_i(d,k)=\binom{k-1}{i}+\binom{k-2}{i-1}$, so
$\phi_{i-1}(d-1,k-1)-\alpha_i(d,k)=\binom{k-2}{i-2}$.

Thus we get the coefficient $A_j$ for $j\ge m+1$,
$$ A_j=(-1)^{d-1-j}
    \sum_{i=0}^{d-1-j}(-1)^i\left[\binom{d-1-i}{j}\binom{k-1}{i}+
    \binom{d-2-i}{j-1}\binom{k-2}{i-1}\right].$$

An induction proof shows that 
$A_j=\binom{k-j-1}{d-j-1}$ for $j\ge m+1$.

Now let $i\le m$; then $d-i\ge m+1$, so
\begin{eqnarray*}
h_i(P^{d,k,n})-h_i(P^{d,k,k})&=&h_{d-i}(P^{d,k,n})-h_{d-i}(P^{d,k,k})\\
&=&(n-k)\binom{k-d+i-1}{i-1}.
\end{eqnarray*}
Since $h_i(P^{d,k,k})=\binom{k-d+i}{i}$ for $i\le m$,
$$h_i(P^{d,k,n})=\binom{k-d+i}{i}+(n-k)\binom{k-d+i-1}{i-1}.$$

\vspace*{-24pt}

\end{proof}

Among all $d$-polytopes with $n+1$ vertices, the minimum $h$-vector 
has $h_i=n-d+1$ for all $i$, $1\le i\le d-1$; this 
is the $h$-vector of the multiplex.
The maximum $h$-vector (for $d$ odd) 
has $h_i=\binom{n-d+i}{i}$ for $1\le i\le (d-1)/2$; this
is the $h$-vector of the cyclic polytope.
The ordinary polytopes provide a nice distribution of $h$-vectors between
these two extremes.
Note that the $g$-vector satisfies the relation 
$g_i(P^{d,k,n})=h_i(P^{d,k-1,n-1})$.

\section{Multiplicial Polytopes with a Large Facet}

The facets of ordinary $d$-polytopes are ``small'': each has at most
$2d-2$ vertices.
Bisztriczky raised the question of whether there exist order-multiplicial
polytopes with a large proportion of the vertices on a single facet.
The corresponding question for weakly multiplicial polytopes is trivial,
since the pyramid over a multiplex is a weakly multiplicial polytope
(but is not generally an order-multiplicial polytope).
We give here a construction of order-multiplicial 4-polytopes having
two-thirds of their vertices on one facet.
\begin{theorem}
For every integer $q\ge 5$, there is an order-multiplicial 4-polytope with
$q$ vertices having a facet containing $\lceil (2q+2)/3 \rceil$ vertices.
\end{theorem}

\begin{proof}
First note that the result is obtained for $q=5$ and $q=6$ by
multiplexes (the simplex, and the double pyramid over a square, respectively).  
So assume $q\ge 7$.
We show in detail the construction for $q\equiv 1\pmod{3}$, and then
give the modifications for the other congruence classes.
So let $q\equiv 1\pmod{3}$, $q\ge 7$.
Let $n=(2q+1)/3$. 
Let $M^{3,n}$ be a 3-multiplex with vertices labeled
by the integers \mbox{$\{0,1, 2, \ldots, q-1\}\setminus \{3j: 1\le j\le (q-2)/3\}$},
with the standard ($<$) order.
Let $Q_0$ be a pyramid over the 3-multiplex with apex labeled 3, 
$Q_0=P(M^{3,n})$.
Now $Q_0$ is (weakly) multiplicial, but not order-multiplicial.  
We perform some subdivisions on $Q_0$ to achieve order-multipliciality.
The polytope $Q_0$ has three kinds of facets: the multiplex $M^{3,n}$ itself,
simplices, and pyramids over quadrilateral faces of $M^{3,n}$.
The multiplex and the simplices are multiplexes in the induced order;
the pyramids in general are not.
Note that the ``first'' two pyramidal facets are multiplexes in the
induced order; their vertex sets are $\{0,1,3, 4, 5\}$ and $\{1, 2, 3, 5, 7\}$
(or, rather, $\{1, 2, 3, 5, 6\}$, if $q=7$).

We pair up the remaining pyramidal facets and add a new vertex beyond each
pair.
For each $i$, $0\le i\le (q-10)/3$, let $F_i$ be the convex hull of 
$\{2+3i, 4+3i, 7+3i, 8+3i, 3\}$, and let $G_i$ be the convex hull of 
$\{4+3i, 5+3i, 8+3i, 10+3i, 3\}$ (or, rather, 
$G_i=\{4+3i, 5+3i, 8+3i, 9+3i, 3\}$, if $i=(q-10)/3$).
Note that these two facets share the 2-face with vertex set $\{4+3i, 8+3i, 3\}$.
Choose a point beyond $F_i$ and $G_i$ and beneath all other facets of $Q_0$, 
and label it $6+3i$.
Let $Q_{i+1}$ be the convex hull of $Q_i$ with the new vertex $6+3i$.
By \cite{Altshuler-Shemer}, the facets of $Q_{i+1}$ consist of the facets of 
$Q_i$, except $F_i$ and $G_i$, and eight additional facets, with vertex 
sets
$\{2+3i, 4+3i, 6+3i, 3\}$,
$\{2+3i, 6+3i, 7+3i, 3\}$,
$\{6+3i, 7+3i, 8+3i, 3\}$,
$\{2+3i, 4+3i, 6+3i, 7+3i, 8+3i\}$,
$\{4+3i, 5+3i, 6+3i, 3\}$,
$\{5+3i, 6+3i, 10+3i, 3\}$,
$\{6+3i, 8+3i, 10+3i, 3\}$, and
$\{4+3i, 5+3i, 6+3i, 8+3i, 10+3i\}$.
(Here, $10+3i$ should be replaced by $9+3i$ in each set, if $i=(q-10)/3$.)
After performing this operation for all $i$, $0\le i\le (q-10)/3$, all
the pyramidal facets of $Q_0$ are replaced by collections of simplices
(necessarily multiplexes), and by pyramidal facets that
are now multiplexes with the induced order.
Thus the resulting polytope, $P=Q_{(q-10)/3+1}$ is an order-multiplicial
polytope.
The polytope $Q_0$ has $n+2$ vertices, and $(q-10)/3+1$ vertices have
been added.  
The polytope $P$ has exactly $q$ vertices.
Of those $q$ vertices, $n+1=(2q+4)/3=\lceil (2q+2)/3 \rceil$ are vertices of
the original base facet $M^{3,n}$.

For $q\equiv 2\pmod{3}$ ($q\ge 8$), let $n=(2q-1)/3$.
We proceed as before with $Q_0=P(M^{3,n})$,
adding vertices labeled $6+3i$ for $0\le i\le (q-11)/3$
for each facet pair $F_i$, $G_i$.
Finally we choose a point beyond the ``last'' facet $F$, which is a simplex
with  vertex set
$\{q-4, q-3, q-1, 3\}$, and beneath all other facets, and label it $q-2$.
Let $P$ be the convex hull of $Q_{(q-11)/3+1}$ with the new vertex $q-2$.
The facets of $P$ consist of the facets of $Q_{(q-11)/3+1}$, except $F$,
and four additional facets, all of which are simplices.
The resulting polytope $P$ is an order-multiplicial polytope with exactly
$q$ vertices.
Of those $q$ vertices, $n+1=(2q+2)/3=\lceil (2q+2)/3 \rceil$ are vertices of
the original base facet $M^{3,n}$.

If $q\equiv 0\pmod{3}$ ($q\ge 9$), let $n=2q/3$.
The polytope $Q_0=P(M^{3,n})$ has an odd number of pyramidal facets.
We proceed as before, adding vertices labeled $6+3i$ for $0\le i\le (q-12)/3$
for each facet pair $F_i$, $G_i$.
The last pyramidal facet $F$ has vertex set $\{q-7, q-5, q-2, q-1, 3\}$.
Choose a point beyond $F$ and beneath all other facets of $Q_{(q-12)/3+1}$,
and label it $q-3$.
Let $P$ be the convex hull of $Q_{(q-12)/3+1}$ with the new vertex $q-3$.
The facets of $P$ consist of the facets of $Q_{(q-12)/3+1}$, except $F$,
and five additional facets, with vertex sets
$\{q-7, q-5, q-3, 3\}$,
$\{q-7, q-3, q-2, 3\}$,
$\{q-5, q-3, q-1, 3\}$,
$\{q-3, q-2, q-1, 3\}$, and
$\{q-7, q-5, q-3, q-2, q-1\}$.
The resulting polytope $P$ is an order-multiplicial polytope with exactly
$q$ vertices.
Of those $q$ vertices, $n+1 = (2q+3)/3=\lceil (2q+2)/3 \rceil$ are vertices
of the original base facet $M^{3,n}$.

\end{proof}

For $q\ge 7$, $q\not\equiv 0\pmod{3}$, the $f$-vector of this multiplicial
4-polytope is
$$f(P)=\left(q,\left\lfloor \frac{13q-43}{3}\right\rfloor,6q-26,
          \left\lfloor \frac{8q-34}{3}\right\rfloor\right)$$
and $f_{02}=\lfloor(56q-242)/3\rfloor$.
For $q\ge 9$, $q\equiv 0\pmod{3}$, the $f$-vector is
$$f(P)=\left(q, \frac{13q-45}{3},6q-27,
           \frac{8q-36}{3}\right)$$
and $f_{02}=(56q-252)/3$.
The $h$-vector is $h=(1,q-4,2q-12,q-4,1)$, unless 
$q\equiv 1\pmod{3}$, in which case it is $h=(1,q-4,2q-11,q-4,1)$.

\section{Conclusion}
We hope that multiplicial polytopes can be used to further our
understanding of the set of all flag vectors of polytopes.
Ordinary polytopes seem particularly important because they generalize
cyclic polytopes, which have played a major role in the study
of simplicial polytopes (see \cite{billera-lee,grunbaum,McMullen,ramsan}).
From the viewpoint of flag vectors, the restriction to odd dimensions is
unfortunate.
However, even-dimensional polytopes can be generated by taking vertex
figures of the odd-dimensional ordinary polytopes.  The resulting polytopes
are multiplicial.


\begin{thebibliography}{10}
\bibitem{Altshuler-Shemer}
Amos Altshuler and Ido Shemer.
\newblock Construction theorems for polytopes.
\newblock {\em Israel J. Math.}, {\bf 47} (1984), 99--110. 

\bibitem{Bayer-Bru-Ste}
Margaret M. Bayer, Aaron M. Bruening, and Joshua D. Stewart,
\newblock A combinatorial study of multiplexes and ordinary polytopes,
\newblock {\em Discrete Comput.\ Geom.}, {\bf 27} (2002), no.~1, 49--63.

\bibitem{bay-ehr}
Margaret~M. Bayer and Richard Ehrenborg,
\newblock The toric $h$-vectors of partially ordered sets, 
\newblock {\em Trans.\ Amer.\ Math.\ Soc.} {\bf 352} (2000), 4515--4531.
    
\bibitem{billera-lee}   
Louis~J. Billera and Carl~W. Lee,
\newblock A proof of the sufficiency of {M}c{M}ullen's conditions for
  $f$-vectors of simplicial polytopes,
\newblock {\em J.\ Combin.\ Theory Ser.\ A} {\bf 31} (1981), 237--255.

\bibitem{Bisz-mult} 
T.~Bisztriczky.
\newblock On a class of generalized simplices.
\newblock {\em Mathematika}, {\bf 43} (1996), 274--285.

\bibitem{Bisz}
T.~Bisztriczky.
\newblock Ordinary $(2m+1)$-polytopes.
\newblock {\em Israel J. Math.}, 102 (1997), 101--123.

\bibitem{Dinh}
Thi~Ngoc Dinh.
\newblock {\em Ordinary Polytopes}.
\newblock PhD thesis, The University of Calgary, 1999.

\bibitem{grunbaum}
Branko Gr\"{u}nbaum.
\newblock {\em Convex Polytopes, Second Edition},
\newblock Springer-Verlag, New York, 2003.

\bibitem{McMullen}
P.~McMullen and G. C. Shephard, 
\newblock {\em Convex polytopes and the upper bound conjecture}, 
\newblock Cambridge Univ. Press, London, 1971.

\bibitem{ramsan}
J{\"o}rg Rambau and Francisco Santos.
\newblock The generalized {B}aues problem for cyclic polytopes. {I}.
\newblock European J. Combin., {\bf 21} (2000), 65--83.

\bibitem{stanley}
Richard~P.\ Stanley,
\newblock The number of faces of simplicial convex polytopes,
\newblock {\em Adv.\ Math.} {\bf 35} (1980), 236--238.

\bibitem{sta87}
Richard~P. Stanley.
\newblock {G}eneralized ${H}$-vectors, intersection cohomology of toric
  varieties, and related results.
\newblock In {\em {C}ommutative algebra and combinatorics ({K}yoto, 1985)},
  volume~11 of {\em Adv. Stud. Pure Math.}, pages 187--213, Amsterdam-New York,
  1987. North-Holland.


\end{thebibliography}
\end{document}